\documentclass[final,date,omega,math,twocoltoc]{diplomega} 

  \usepackage{amsmath}
  \usepackage{amsfonts}
  \usepackage{amssymb}

  \usepackage{url}
  \usepackage{nicefrac}

  \usepackage[final]{epsfig}
  \usepackage[final]{graphicx}
  \usepackage{tabularx}
  \usepackage[small]{caption}
  \usepackage[all,poly,dvips]{xy}
  \CompileMatrices
  \SelectTips{cm}{}

  \usepackage{pifont}

\makeindex
\newindex{not}{ndx}{nnd}{Table of Notation}

\newcommand{\Gromovconjecture}{%
  \begin{conj}
    Let $M$ be an aspherical closed orientable manifold 
    whose simplicial volume vanishes. Then for all~$p \in \N$,
    \[ \ltb[big]{p}{\ucov{M};{\mathcal N}(\pi_1(M))} = 0.\]
 \end{conj}%
}
\newcommand{\ltproptmp}{%
  Let $U$ be a simply connected Riemannian manifold. Then 
  for each $p \in\N$ there are 
  constants~$B_p \in \R_{> 0}$ and a constant~$T \in \R_{>0}$
  satisfying: for any discrete group~$G$ acting
  freely and cocompactly on~$U$ by isometries and all~$p \in \N$, 
  \begin{align*}
    \ltb[big]{p}{U; 
                  {\mathcal N}(G)} 
   = B_p \cdot \vol(G\setminus U).
  \end{align*}
  If $U$ with this $G$-action is of determinant class, then
  \[ \ltt[big]{U; {\mathcal N}(G)} = T \cdot \vol(G\setminus U).\]
}
\newcommand{\ltpropprincunn}{%
   \begin{satzunn}[Proportionality Principle of $L^2$-Invariants]
     \ltproptmp
   \end{satzunn}
}

\newcommand{\svproptmp}{%
  Let $M$ and $N$ be oriented closed connected Riemannian manifolds
  with isometrically isomorphic universal
  Riemannian coverings. Then 
  \[ \frac{\sv M}{\vol (M)} = \frac{\sv N}{\vol (N)}.\]}

\newcommand{\svpropprincunn}{%
  \begin{satzunn}[Proportionality Principle of Simplicial Volume]
    \svproptmp
  \end{satzunn}
}
\newcommand{\hypproptmp}{%
  For $n \in\N$ define $v_n \in \R_{>0}$ as the maximal volume of an
  ideal $n$-simplex in $n$-dimensional hyperbolic space. 
  Then for each hyperbolic oriented closed connected manifold~$M$ of
  dimension~$n$, 
  \[ \frac{\sv M}{\vol(M)} = \frac1{v_n}. \]
}

\newcommand{\hyppropprincunn}{%
  \begin{satzunn}\hypproptmp
  \end{satzunn}
}

\begin{document}
\frontmatter


\begin{titlepage}
  \rule{\linewidth}{0pt} 
  \vfill
  \begin{flushright}
    \Huge{The Proportionality Principle\\
          of Simplicial Volume}

    \vspace*{1em}

    \Large{Clara Strohm}
  \end{flushright}
  \vfill\vfill
\end{titlepage}
\thispagestyle{empty}
\begin{titlepage}
  \newlength{\wolf}
  \settowidth{\wolf}{Prof.\ Dr.\ Wolfgang L\"uck}
  \rule{\linewidth}{0pt} 
  \begin{flushright}
    Westf\"alische Wilhelms-Universit\"at M\"unster \\
    Fachbereich Mathematik und Informatik
  \end{flushright}
  \vfill
  \begin{flushright}
    \Large{Diplomarbeit}

    \vspace*{1em}

    {\Huge The Proportionality Principle\\
           of  Simplicial Volume}

    \vspace*{1em}

    {\Large Clara Strohm}

    \vspace*{.5em}
     
    \normalsize{September 2004}
   \end{flushright}
   \vfill\vfill
   \begin{flushright}
     \begin{tabular}{rr}
       Betreuer :      & Prof.\ Dr.\ Wolfgang L\"uck \\
       Zweitgutachter: & \makebox[\wolf][s]{PD\ Dr.\ Michael Joachim}
     \end{tabular}
   \end{flushright}
\end{titlepage}
\thispagestyle{empty}
\rule{\linewidth}{0pt} 

\vfill

Hiermit versichere ich, da\ss\ ich diese Arbeit selbst\"andig verfa\ss t
und nur die angegebenen Quellen und Hilfsmittel benutzt habe. 

\vspace*{1em}

M\"unster, im September 2004

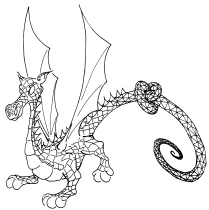
\chapter*{Introduction}

Manifolds, the basic objects of geometry, form the playground for both
topological and smooth structures. Many challenges in modern
mathematics are concerned with the nature of this interaction between
algebraic topology, differential topology and differential
geometry. For example all incarnations of the index theorem, 
the Gau\ss-Bonnet formula, 
and the Mostow
rigidity theorem fall into this category. 

\begin{satzunn}[Mostow Rigidity Theorem]\tindex{Mostow rigidity theorem}%
  Each homotopy equivalence between 
  oriented closed connected hyperbolic manifolds of dimension at least~3
  is homotopic to an isometry. 
\end{satzunn}

This theorem was first proved by Mostow \cite{mostow}.  
But there is also a proof by Gromov \cite{munkholm}, for which he
designed the simplicial volume (therefore also known as ``Gromov
norm'').

The simplicial volume is a homotopy invariant of oriented closed
connected manifolds defined in terms of the singular chain complex
with real coefficients. This invariant consists of a 
nonnegative real number measuring the efficiency of representing
the fundamental class using singular simplices: 
If $M$ is an oriented closed connected manifold, then
the \emph{simplicial volume} of~$M$, denoted by~$\sv M$, is the infimum of
all $\ell^1$-norms of singular cycles with real coefficients representing the
fundamental class.
Since the fundamental class can be viewed as a
generalised triangulation of the manifold, the simplicial volume can
also be interpreted as a measure for the complexity of a manifold. 

But the simplicial volume also has another important interpretation:
it is a homotopy invariant approximation of the Riemannian volume. For
example, 
Gromov's proof of the Mostow rigidity theorem is based on the
following proportionality principle for hyperbolic manifolds 
\cite[Section~2.2]{gromov}, \cite[Theorem~6.2]{thurston}:

\hyppropprincunn

In particular, the simplicial volume does not always vanish.  However,
in general the simplicial volume is rather hard to compute. In the
most cases where the simplicial volume is known, it is zero.  This,
for example, is the case if the fundamental group of the manifold is
Abelian or, more generally, amenable.

One of the most complete (but hardly digestible) references about the
simplicial volume is Gromov's pioneering paper ``Volume and
Bounded Cohomology'' \cite{gromov}. More recent introductions are
given in the textbooks on hyperbolic geometry by Ratcliffe and
Benedetti, Petronio \cite[\S11.5]{ratcliffe}, \cite[Section~C.3]{bp}.

Applications of the simplicial volume also include knot theory
\cite{murakami}, the existence of certain $S^1$-operations
\cite{yanos1}, and the investigation of other volumes, such as the
Riemannian or the minimal volume \cite{gromov}.  
Furthermore, Gromov conjectured the following relation between the
simplicial volume and $L^2$-invariants \cite[Section 8A]{gromovasym}
-- thereby connecting two distant fields of algebraic topology:

\Gromovconjecture

A rather complete account of $L^2$-invariants is given in the book
\emph{$L^2$-Invariants: Theory and Applications to Geometry and
$K$-theory} \cite{lueck}. The above conjecture is justified by the 
observation that $L^2$-invariants and the simplicial volume show a similar
behaviour. We will discuss this conjecture briefly in
Section~\ref{l2sec}. Among the most striking similarities is the
proportionality principle:

\ltpropprincunn

In particular, if $M$ and $N$ are oriented closed connected Riemannian
manifolds with isometrically isomorphic universal Riemannian
coverings, then for all $p \in \N$,
\[ \frac{\ltb[big]{p}{\ucov M; {\mathcal N}(\pi_1(M))}}{\vol(M)}
 = \frac{\ltb[big]{p}{\ucov N; {\mathcal N}(\pi_1(N))}}{\vol(N)}
\]
and
\[ \frac{\ltt[big]{\ucov M; {\mathcal N}(\pi_1(M))}}{\vol(M)}
 = \frac{\ltt[big]{\ucov N; {\mathcal N}(\pi_1(N))}}{\vol(N)}.
\]

\svpropprincunn

These proportionality principles are also examples for an interesting
link between topological and differential structures. A similar
proportionality principle is known for Chern numbers
\cite[Satz~3]{hirzebruch}. 

A proof in the $L^2$-case is given in L\"uck's book
\cite[Theorem~3.183]{lueck}.  For the simplicial volume, sketches of
(different) proofs were given by Gromov \cite[Section~2.3]{gromov} and
Thurston \cite[page~6.9]{thurston}.  However, there is no complete
proof of the proportionality principle of simplicial volume in the
literature. It is the aim of this diploma thesis to provide such a
proof, based on Thurston's approach.

This text is organised as follows:

In the first chapter the simplicial volume is introduced and some of
its properties are collected.

Chapter~2 is devoted to the study of bounded cohomology, a quite
peculiar variant of singular cohomology. Its application as a basic
tool for the analysis of the simplicial volume was discovered by
Gromov \cite{gromov}.  However, we will mostly refer to the more
elegant approach due to Ivanov, based on classical homological algebra
\cite{ivanov}, \cite{monod}.

Thurston's proof of the proportionality principle relies on the
computation of the simplicial volume by means of a new homology theory
called measure homology \cite[pages~6.6--6.9]{thurston}.  Instead of
linear combinations of singular simplices, measures on the set of
(smooth) singular simplices are considered as chains, thereby
permitting averaging constructions such as Thurston's smearing.  But
Thurston's exposition lacks a proof of the fact that measure homology
and singular homology give rise to the same simplicial volume.  We
will define measure homology in Chapter~3 and list its basic algebraic
properties as they can be found in the papers of Zastrow and Hansen 
\cite{zastrow}, \cite{hansen}, in
a slightly different setting. 

Chapter~4 constitutes the central contribution of this diploma
thesis. It will be shown that measure homology indeed can be used to
compute the simplicial volume, i.e., that singular homology with real
coefficients is \emph{isometrically} isomorphic to measure
homology. The proof is based on techniques from bounded cohomology, as
explained in Chapter~2. It will also be discussed why Bowen's argument
\cite{bowen} is not correct (Subsection~\ref{bowen}).

The last chapter finally contains the proof of the proportionality
principle -- using Thurston's smearing technique. This forces us to
study isometry groups of Riemannian manifolds.  En passant, it will be
proved that the compact open topology and the \cone-topology (the
smooth counterpart of the compact open topology) on the isometry
groups coincide. Unfortunately, I was not able to show this without
making use of the theory of standard Borel spaces.  The chapter closes
with a proof of the fact that the simplicial volume of a flat
orientable closed connected manifold vanishes. The proof is based on
the proportionality principle, instead of referring to Gromov's
sophisticated estimate of the simplicial volume via the minimal volume
\cite[page~220]{gromov}.

I tried to keep the prerequisites for the understanding of this diploma
thesis as small as possible. So only a basic familiarity with singular
(co)homology, covering theory, measure theory, and Riemannian geometry
is required. A variety of textbooks is available providing such a
background: e.g., \cite{lee}, \cite{leesmooth}, \cite{massey}, \cite{bredon},
\cite{elstrodt}. But the detailed exposition of the subject comes at a
price -- the length of the text, for which I apologise.

The symbols $\Z$, $\N$, $\Q$, $\R$, $\C$ stand for the set of integers,
the set of nonnegative integers, the set of rational numbers, the set
of real numbers, and the set of complex numbers respectively. The one
point space is written as~$\opsp$. For simplicity, we always assume
that manifolds have dimension at least~1. The (co)homology class
represented by an element~$x$ is denoted by~$[x]$. 
All sums ``$\sum$'' are implicitly finite, unless stated otherwise. 
All notation is
collected in the table of notation (between bibliography and index). 

\snindex{$a}{$\consum$}{connected sum}%
\snindex{$b}{$\opsp$}{one point space}%
\snindex{$c}{$\lor$}{wedge}%
\snindex{ker}{$\ker$}{kernel}%
\snindex{im}{$\im$}{image}%
\snindex{$s}{$[x]$}{(co)homology class represented by~$x$}%
\snindex{vol}{$\vol(M)$}{volume of~$M$}%
\snindex{pi1}{$\pi_1$}{fundamental group}%
\snindex{Tf}{$Tf$}{differential of~$f$}%
\snindex{volM}{$\vol_M$}{volume form of~$M$}%
\snindex{C}{$\C$}{complex numbers}%
\snindex{N}{$\N$}{nonnegative integers}%
\snindex{Q}{$\Q$}{rational numbers}%
\snindex{R}{$\R$}{real numbers}%
\snindex{Z}{$\Z$}{integers}%
\snindex{colim}{$\dirlim$}{colimit/direct limit}%
\vspace*{-\baselineskip}
\section*{Acknowledgements}
\enlargethispage{\baselineskip}

First of all I would like to thank Professor Wolfgang L\"uck for
giving me the opportunity to write this diploma thesis, for the 
encouragement, and for the freedom I enjoyed during this project.

Roman Sauer and Marco Schmidt infected me with their fascination for
the simplicial volume and bounded cohomology. I am very grateful for
all their support and guidance.  Marco Varisco was always open for
discussions on mathematical principles and notation, and for making
XYZ-theory accessible for pedestrians like me. Mille grazie!

Professor Volker Puppe and his concise topology lectures and
seminars led me through my two years at the University of 
Konstanz. I would like to thank him  for all his support. 

Furthermore, I would like to thank Stefanie Helker for studying English and 
Benedikt Plitt, 
Julia Weber (Wohlf\"uhlb\"uro!), and 
Michael Weiermann 
for hunting bugs and for their valuable suggestions.

My parents always supported me, even though I did not manage to 
become acquainted with philosophy -- but at least I am no 
\foreignlanguage{polutonikogreek}{>agewm'etrhtoc} 
anymore. Dankesch\"on!

All my gratitude is devoted to Andres L\"oh, not only for his 
\TeX nical ekspertise, but also for his patience and all the things 
language is not powerful enough to express -- except for
Pl\"uschtiere: Wuiisch? Zonder pf!

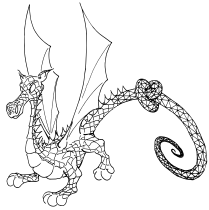
\begingroup
\makeatletter
\let\@mkboth=\@gobbletwo
\makeatother
\let\cleardoublepage=\clearpage
\tableofcontents
\endgroup

%
\mainmatter

\chapter{Simplicial Volume}

\begin{chapterabstract} 
 The simplicial volume of an oriented closed connected manifold is a
 homotopy invariant -- defined in terms of the singular chain complex -- 
 which measures the efficiency of representing the fundamental class by
 singular chains (with real coefficients). Its first appearance is in
 Gromov's famous proof of the Mostow rigidity theorem
 \cite{munkholm}. 

 Despite of being a topological invariant the simplicial volume contains
 interesting information about the possible differential structures on
 a manifold -- in form of the Riemannian volume. Examples for such
 relations are Gromov's estimate of the minimal volume and the
 proportionality principle, which is the topic of this thesis.  
 
 Moreover, there seems to be a deep connection to  
 $L^2$-invariants \cite[Section~8A]{gromovasym}:
 
 \Gromovconjecture
  
 In this chapter the simplicial volume is introduced
 (Section~\ref{svdefsec}) and some of its properties are collected
 (Section~\ref{svpropsec}).  We then give a short survey about
 indications supporting the above conjecture in Section~\ref{l2sec}.
 In the last section we briefly discuss some generalisations of the
 simplicial volume.
\end{chapterabstract}

\InputIfFileExists{simvol}{}


\chapter{Bounded Cohomology}

\begin{chapterabstract}\label{chboco}
 Bounded cohomology is the functional analytic twin of singular
 cohomology -- it is constructed similarly to singular cohomology from
 the singular chain complex, using the topological dual space instead of the
 algebraic dual space \cite{gromov}. The corresponding duality on the
 level of (co)homology is explored in Section~\ref{dualitysection}.
 Using this duality, it is possible to calculate the simplicial volume
 via bounded cohomology.
 
 It is astonishing that the rather small difference between the
 definitions of singular and bounded cohomology lead to
 completely different characters of both theories. 

 Apart from the above geometric definition, there is also a notion of
 bounded cohomology of groups due to Ivanov based on techniques from
 homological algebra \cite{ivanov}, which we will sketch in
 Section~\ref{bcgroupssection}. This algebraic approach is the source
 of the strength of bounded cohomology.

 In Section~\ref{mtsec}, a key feature of bounded cohomology is
 discussed: the bounded cohomology of a topological space coincides
 with the bounded cohomology of its fundamental group.  Unfortunately,
 unlike singular cohomology bounded cohomology fails to satisfy the
 excision axiom. Nevertheless, in many cases bounded cohomology can be
 calculated directly from special resolutions. For example, it can be
 shown that bounded cohomology ignores amenable groups, implying that
 the bounded cohomology of spaces with amenable fundamental group
 vanishes.
 
 In the last section, we will explain a special resolution for
 calculating the bounded cohomology of a group, which plays a crucial
 r\^ole in Section~\ref{isometricisom}. 
\end{chapterabstract}

\InputIfFileExists{boco}{}

\addtocontents{toc}{\protect\end{multicols}}
\addtocontents{toc}{\protect\newpage}
\addtocontents{toc}{\protect\rule{\linewidth}{0pt}}
\addtocontents{toc}{\protect\vspace*{12em}}
\addtocontents{toc}{\protect\begin{multicols}{2}}
\chapter{Measure Homology}\label{chmh}

\begin{chapterabstract}
  Measure homology is a generalisation of singular homology in the
  following way: singular chains with real coefficients are viewed
  as signed measures on the space of singular simplices whose mass is
  concentrated in finitely many points. In the measure homology chain 
  complex also more complicated measures are allowed.

  Measure homology was introduced in Thurston's lecture notes
  in~1979 \cite[page~6.6]{thurston}. He already claimed that
  measure homology and singular homology should coincide. 
  More extensive accounts are the papers of Zastrow \cite{zastrow} and
  Hansen \cite{hansen} and the book of Ratcliffe \cite{ratcliffe}.
  
  The motivation for the introduction of measure homology originates
  from the fact that measure homology can be used to calculate the
  simplicial volume, hence giving more room for constructions such as
  smearing (cf.~Chapter~\ref{chpropprinc}). This smearing construction
  will be the key to the proof of the proportionality principle.

  In Section~\ref{prelude}, the technical background needed for the
  definition of measure homology is presented -- in particular, a
  feasible topology on the set of smooth singular simplices will be
  introduced. In Section~\ref{mhdefsec} measure homology is
  defined. The basic properties of measure homology are listed in
  Section~\ref{mhprop}, and there is also given a proof for the
  compatibility of measure homology with colimits.  The next chapter
  is devoted to the proof that measure homology and singular homology
  are isometrically isomorphic.
\end{chapterabstract}

\InputIfFileExists{meashom}{}


\chapter{The Proportionality Principle}\label{chpropprinc}

\begin{chapterabstract}
  The proportionality principle of simplicial volume reveals a
  fascinating connection between the simplicial volume and the
  Riemannian volume: simplicial volume and Riemannian volume are
  proportional if we restrict our attention to manifolds which share
  the same universal Riemannian covering space. 
  Similar proportionality principles also occur in the setting of
  $L^2$-invariants (Theorem~\ref{l2pp}) and in the setting of
  characteristic numbers of complex manifolds
  \cite[Satz~3]{hirzebruch}. 

  Both Thurston and Gromov sketched (dual) proofs for the
  proportionality principle of simplicial volume \cite[page~6.9]{thurston},
  \cite[Section~2.3]{gromov}. Thurston's idea is to use measure
  homology to compute the simplicial volume and to take advantage of
  the larger chain complex for the so-called smearing construction.
  In this chapter, a detailed version of Thurston's proof will be
  given, based on the results of the previous chapters.

  Gromov's proof uses arguments from bounded cohomology and also
  depends on some averaging operation via Haar measures. Unfortunately his
  exposition is not very explicit about certain measurability issues. 
 
  In Section~\ref{stpp}, the proportionality principle is discussed and
  the strategy of proof is presented. Sections \ref{isometrygroups}
  and \ref{haarmeasure} provide the necessary tools for Thurston's
  smearing technique (i.e., the study of isometry groups and
  appropriate measures on them).  In Section~\ref{smearing}, the
  smearing of smooth singular chains is developed. The actual proof
  of the proportionality principle is given in
  Section~\ref{pproof}. The chapter ends with some applications in
  Section~\ref{ppapp}.
\end{chapterabstract}

\InputIfFileExists{propprinc}{}


\printindex[not]

\cleardoublepage
\printindex

\end{document}